\documentclass[a4paper,12pt]{amsart}
\usepackage{amsfonts}
\usepackage{amsmath}
\usepackage{amssymb}
\usepackage[a4paper]{geometry}
\usepackage{mathrsfs}
\usepackage{hyperref}
\theoremstyle{plain}
\newtheorem{thm}{Theorem}[section]
\newtheorem{prop}[thm]{Proposition}
\newtheorem{cor}[thm]{Corollary}
\newtheorem{lem}[thm]{Lemma}
\theoremstyle{definition}
\newtheorem{defn}[thm]{Definition}

\newtheorem{ex}[thm]{Example}
\def\G{G}
\def\U{\mathrm U}
\def\C{\mathbb C}
\DeclareMathOperator{\spann}{span}
\DeclareMathOperator{\trace}{Tr}
\def\d{\mathrm d}
\DeclareMathOperator{\singsupp}{sing\,supp}
\def\diff{\mathrm{diff}}
\def\Diff{\mathrm{Diff}}
\DeclareMathOperator{\rank}{rank}
\def\Gh{{\widehat{\G}}}
\def\SU2{{{\rm SU}(2)}}
\def\irm{{\rm i}}
\def\erm{{\ {\rm e}}}
\def\drm{{\ {\rm d}}}
\def\R{{\mathbb R}}
\def\D{{\mathbb D}}
\def\Darr{{\overrightarrow{\D}}}
\def\N{{\mathbb N}}
\def\Z{{\mathbb Z}}
\def\C{{\mathbb C}}
\def\S3{{{\mathbb S}^3}}
\def\Rn{{{\mathbb R}^n}}
\def\Tn{{{\mathbb T}^n}}
\def\Zn{{{\mathbb Z}^n}}

\def\Lap{{\mathcal L}}

\def\Op{{{\rm Op}}}

\def\p#1{{\left({#1}\right)}}

\def\br#1{{\left[{#1}\right]}}
\def\jp#1{{\left\langle{#1}\right\rangle}}
\def\n#1{{\left\|{#1}\right\|}}

\def\Smrd{{\mathscr S^{m}_{\rho,\delta}(G)}}
\begin{document}
\title[H\"ormander class $\Psi^m(\G)$ on compact Lie groups]
{H\"ormander class of pseudo-differential operators
on compact Lie groups and global hypoellipticity}
\author[Michael Ruzhansky]{Michael Ruzhansky}
\address{
  Michael Ruzhansky:
  \endgraf
  Department of Mathematics
  \endgraf
  Imperial College London
  \endgraf
  180 Queen's Gate, London SW7 2AZ 
  \endgraf
  United Kingdom
  \endgraf
  {\it E-mail address} {\rm m.ruzhansky@imperial.ac.uk}
  }
\author[Ville Turunen]{Ville Turunen}
\address{
  Ville Turunen:
  \endgraf
   Helsinki University of Technology
  \endgraf
  Institute of Mathematics
  \endgraf
   P.O. Box 1100
  \endgraf
   FI-02015 HUT
  \endgraf
  Finland
  \endgraf
  {\it E-mail address} {\rm ville.turunen@hut.fi}
  }  
\author[Jens Wirth]{Jens Wirth}
\address{
  Jens Wirth:
  \endgraf
  Department of Mathematics
  \endgraf
  Imperial College London
  \endgraf
  180 Queen's Gate, London SW7 2AZ 
  \endgraf
  United Kingdom
  \endgraf
  {\it E-mail address} {\rm j.wirth@imperial.ac.uk}
  }
\thanks{The first
 author was supported by the EPSRC
 Leadership Fellowship EP/G007233/1. The third author
 was supported by the EPSRC grant EP/E062873/1.
 }
\date{\today}


\begin{abstract}
In this paper we give several global characterisations 
of the H\"or\-mander class $\Psi^m(\G)$ of
pseudo-differential operators on compact Lie groups.
The result is applied to give criteria for the ellipticity
and the global hypoellipticity of pseudo-differential operators 
in terms of their matrix-valued full symbols.
Several examples of the first and second order globally
hypoelliptic differential operators are given. Where the
global hypoelliptiticy fails, one can construct explicit
examples based on the analysis of the global symbols.
\end{abstract}

\maketitle


\section{Introduction}

\medskip

\subsection{}
For a compact manifold $M$ we denote by
$\Psi^m(M)$ the set of H\"ormander's pseudo-differential 
operators on $M$, i.e., the class of operators which
in all local coordinate charts give operators in
$\Psi^m(\Rn)$. Operators in $\Psi^m(\Rn)$ are characterised
by the symbols satisfying
$$
|\partial_\xi^\alpha\partial_x^\beta a(x,\xi)|\leq
C(1+|\xi|)^{m-|\alpha|}
$$
for all multi-indices $\alpha, \beta$ and all
$x,\xi\in\Rn$. An operator in $\Psi^m(M)$ is called elliptic
if all of its localisations are locally elliptic. 
The principal symbol of a pseudo-differential operator
on $M$ can be invariantly defined as function on the
cotangent bundle $T^*M$, but it is not possible to control
lower order terms in the same way.
If one fixes a connection, however, it is possible to make
sense of a full symbol, see Widom \cite{Wi80}.

However, on a Lie group $G$ it is a very natural idea
to use the Fourier analysis and the global
Fourier transform on the group to analyse
the behaviour of operators which are originally defined
by their localisations. This allows one to make a full use
of the representation theory and of many results available
from the harmonic analysis on groups. In \cite{RT-groups}
and \cite{RT-book}, the first two authors carried out a
comprehensive non-commutative
analysis of an analogue of the Kohn--Nirenberg
quantisation on $\Rn$, in terms of the representation
theory of the group. This yields a full matrix-valued symbol
defined on the product $G\times\Gh$, which can be viewed
as a non-commutative version of the phase space. The calculus and 
other properties of this quantisation resemble those 
well-known in the local theory on $\Rn$. However, due to
the non-commutativity of the group in general, the full symbol is
matrix-valued, with the dimension of the matrix symbol
$a(x,\xi)$ at $x\in G$ and $[\xi]\in\Gh$ equal to the
dimension of the representation $\xi$.
This construction resembles that of Taylor \cite{Tay-noncomm}
but is carried out entirely in terms of the group $G$ without
referring to its Lie algebra and the Euclidean
pseudo-differential classes there. 
In \cite{Beals} (see also \cite{Taylor97}), pseudo-differential 
operators have been characterised by their commutator properties
with the vector fields. In \cite{RT-book}, we gave a 
different version of this characterisation relying on the
commutator properties in Sobolev spaces, which led 
in \cite{RT-groups} (see also \cite{RT-book}) to
characterisations of symbols of operators in $\Psi^m(G)$.
However, the results there still rely on 
commutator properties of operators and some of them are taken as
assumptions. One aim of this paper is to eliminate such
assumptions from the characterisation, and to give its
different versions relying on different choices of
difference operators on $\Gh$. This will be given in
Theorem \ref{thm:main}.

In Theorem \ref{THM:SU2thm} we give another simple
description of H\"ormander's classes on the group $\SU2$
and on the 3-sphere $\S3$. This is possible due to 
the explicit analysis carried out in \cite{RT-groups}
and the explicit knowledge of the representation of $\SU2$.
We note that this approach works globally on the whole
sphere, e.g. compared to the analysis of Sherman
\cite{She75} working only on the hemisphere
(see also \cite{She90}).
Using the classification of representations in
some cases (see, e.g., \cite{SV}) it is possible to draw
conclusions of this type on other groups as well.
However, this falls outside the scope of this paper.

We give two applications of the characterisation provided by
Theorem \ref{thm:main}. First, we
characterise elliptic operators in $\Psi^m(G)$ by
their matrix-valued symbols. 
Similar to the toroidal case (see \cite{Agranovich2})
this can be applied further to spectral problems.
Second, we give a
sufficient condition for the global hypoellipticity of
pseudo-differential operators. Since the hypoellipticity
depends on the lower order terms of an operator, a 
knowledge of the full symbol becomes crucial. 
Here, we note that while classes of hypoelliptic symbols
can be invariantly defined on manifolds by localisation
(see, e.g., Shubin \cite{Shubin}), the lower order terms
of symbols can not. From this point of view conditions of Theorem
\ref{THM:RTW-hypoellipticity} appear natural as they
refer to the full symbol defined globally on $G\times\Gh$.

The global hypoellipticity of vector fields and second order
differential
operators on the torus has been extensively studied in the
literature, see e.g. \cite{BZS}, \cite{HP}, and references
therein. In this paper we study the global hypoellipticity
problem for pseudo-differential operators on general compact
Lie groups, and discuss a number of explicit examples on the group
$\S3\cong{\rm SU}(2)$. 
We say that an operator $A$ is globally hypoelliptic on $G$
if $Au=f$ and $f\in C^\infty(G)$ imply $u\in C^\infty(G)$.
Thus, we give examples of first
and second order differential operators which are not locally
hypoelliptic and are not
covered by H\"ormander's sum of the squares theorem but
which can be seen to be globally hypoelliptic by the
techniques described in this paper. 
For example, if $X$ is a left-invariant vector field on
$\S3$ and 
$\partial_X$ is the derivative with respect to $X$,
then the operator 
$\partial_X+c$ is globally hypoelliptic
if and only if $\mathrm i c \not\in\frac12\mathbb Z$.
We note that this operator is not locally hypoelliptic,
but the knowledge of the global symbol allows us to draw
conclusions about its global properties. It is interesting to
observe that the number theory plays a role in the global
properties, in some way similar to the appearance of the
Liouville numbers in the global properties of vector fields
on the torus. For example, the failure of the
D'Alembertian $\mathrm D_3^2-\mathrm D_1^2-\mathrm D_2^2$
to be globally hypoelliptic can be related to properties of
the so-called triangular numbers, and an explicit counterexample
to the global hypoellipticity can be consequently constructed
based on these numbers.
We also
look at the cases of the sub-Laplacian and other operators, 
which are covered
by H\"ormander's theorem, but for which we construct
an explicit inverse which turns out to be
a pseudo-differential operator
with the global matrix-valued symbol of type $(\frac12,0)$.
Such examples can also be extended to variable coefficient
versions.

Another application of Theorem \ref{thm:main} to obtain the sharp
G{\aa}rding inequality on Lie groups will appear
elsewhere.

\subsection{} We will now introduce some notation. 
Let throughout this paper $\G$ be a compact Lie group 
of (real) dimension $n=\dim\G$ with the unit element
$e$, and denote by $\widehat\G$ the set 
of all equivalence classes of 
continuous
irreducible unitary representations of $\G$. 
For necessary details on Lie groups and their representations 
we refer to \cite{RT-book}, but recall some basic 
facts for the convenience of the reader and
in order to fix the notation. 
Each $[\xi]\in \widehat\G$ corresponds to a homomorphism 
$\xi : \G \to \U(d_\xi)$ with $\xi(xy)=\xi(x)\xi(y)$ 
satisfying the irreducibility condition $\C^{d_\xi} = 
\spann\{ \xi(x) v : x\in\G\}$ for any given 
$v\in \C^{d_\xi}\setminus\{0\}$. The number $d_\xi$ 
is referred to as the dimension of the representation $\xi$.

We always understand the group $\G$ as 
a manifold endowed with the (normalised) bi-invariant 
Riemannian structure. Of major interest for us is the 
following version of the Peter-Weyl theorem defining 
the group Fourier transform
$\mathcal F : L^2(\G) \to \ell^2(\widehat\G)$.
\begin{thm}\label{thmPW}
The space $L^2(\G)$ decomposes as the
orthogonal direct sum of bi-invariant 
subspaces parameterised by $\widehat\G$, 
$$
   L^2(\G) = \bigoplus_{[\xi]\in\widehat\G} 
   \mathcal H^\xi,\qquad \mathcal H^\xi = 
   \{ x\mapsto \trace (A\xi(x)) : A\in \C^{d_\xi\times d_\xi}\},
$$
the decomposition given by the Fourier series
$$
  f(x) = \sum_{[\xi]\in\widehat\G} d_\xi \trace(\xi(x) 
  \widehat f(\xi)),\qquad \widehat f(\xi) = 
  \int_\G f(x) \xi^*(x) \d x.
$$
Furthermore, the following Parseval identity is valid,
$$
  \|f\|_{L^2(\G)}^2 = \sum_{[\xi]\in\widehat\G} 
  d_\xi \| \widehat f(\xi) \|_{HS}^2. 
$$
\end{thm}

The notion of Fourier series extends naturally to 
$C^\infty(\G)$ and to the space of distributions 
$\mathcal D'(\G)$ with convergence in the respective 
topologies. Any operator $A$ on $\G$ mapping 
$C^\infty(\G)$ to $\mathcal D'(\G)$ gives rise to a 
{\em matrix-valued full symbol} 
$\sigma_A(x,\xi)\in \C^{d_\xi\times d_\xi}$
$$
  \sigma_A(x,\xi) = \xi^*(x) (A\xi)(x),
$$
so that 
\begin{equation}\label{EQ:RTW-quant}
  A f(x) = \sum_{[\xi]\in\widehat\G} d_\xi 
  \trace\big(\xi(x) \, \sigma_A(x,\xi) \, \widehat f(\xi)\big)
\end{equation} 
holds as $\mathcal D'$-convergent series. 
For such operators we will also write $A=\Op(\sigma_A)$.
For a rather comprehensive treatment of this quantisation
we refer to \cite{RT-book} and \cite{RT-groups}.
We denote the right-convolution kernel of $A$ by $R_A$, so that
$$
 Af(x)=\int_\G K_A(x,y)\ f(y) \drm y
 =\int_\G f(y)\ R_A(x,y^{-1}x) \drm y,
$$
where $\drm y$ is the normalised
Haar measure on $G$. The symbol $\sigma_A$
and the right-convolution kernel $R_A$ are related by
$
  \sigma_A(x,\xi) = \int_G R_A(x,y)\ \xi^*(y)\ {\rm d}y.
$

The first aim of this paper 
is to characterise H\"ormander's class of 
pseudo-differential operators $\Psi^m(\G)$ 
by the behaviour of their full symbols, with applications
to the ellipticity and global hypoellipticity of operators.

The plan of this note is as follows. In the next section 
we give several equivalent symbolic characterisations of 
pseudo-differential operators by their symbols. These 
characterisations are based on difference operators 
acting on sequences of matrices 
$\widehat\G\ni[\xi]  \mapsto  \sigma(\xi) \in 
\C^{d_\xi\times d_\xi}$ of varying dimension. 
We also give an application to pseudo-differential operators
on the group $\SU2$.
Section~\ref{sec3} contains the proof and some useful 
lemmata on difference operators. 
Section~\ref{sec4} contains the criteria for the ellipticity 
of operators in terms of their full matrix-valued symbols
as well as a finite version of the Leibniz formula for
difference operators associated to the group
representations. Finally, Section~\ref{sec5} contains
criteria for the global hypoellipticity of 
pseudo-differential operators and a number of examples.


\section{Characterisations of H\"ormander's class}\label{sec2}
The following statement is based on difference operators. 
They are defined as follows. 
We say that $Q_\xi$ is a {\em difference operator} 
of order $k$ if it is given by
$$
   Q_\xi \widehat f(\xi) = \widehat{q_Q f}(\xi)
$$
for a function $q=q_Q\in C^\infty(\G)$ vanishing of order $k$ at
the identity $e\in\G$, i.e., $(P_x q_Q)(e)=0$ for 
all left-invariant differential operators 
$P_x\in\Diff^{k-1}(\G)$ of order $k-1$. 
We denote the set of all difference operators of order $k$ as 
$\diff^k(\Gh)$.

\begin{defn}
A collection of $m\ge n$ first order difference 
operators 
$\triangle_1, \ldots, \triangle_m\in\diff^1(\widehat\G)$
is called {\em admissible}, if the corresponding 
functions $q_1, \ldots, q_m\in C^\infty(\G)$ satisfy 
$\d q_j(e)\ne 0$, $j=1,\ldots,m$, and
$\rank(\d q_1(e),\ldots,\d q_m(e))=n$. 
It follows, in particular, that $e$ is an isolated common zero
of the family $\{q_j\}_{j=1}^m$.
An admissible collection is called {\em strongly admissible} if 
$\bigcap_j\{ x\in\G : q_j(x)=0\} = \{e\}$. 
\end{defn} 

For a given admissible selection of difference operators 
on a compact Lie group $G$ we use multi-index notation 
$\triangle_\xi^\alpha = \triangle_1^{\alpha_1}\cdots 
\triangle_m^{\alpha_m}$ and 
$q^\alpha(x) = q_1(x)^{\alpha_1}\cdots q_m(x)^{\alpha_m}$. 
Furthermore, there exist corresponding differential operators 
$\partial_x^{(\alpha)}\in \Diff^{|\alpha|}(\G)$ 
such that Taylor's formula
\begin{equation}\label{EQ:RTW-Taylor-exp}
   f(x) = \sum_{|\alpha|\le N-1} \frac1{\alpha!} \, 
   q^\alpha(x) \, \partial_x^{(\alpha)} f(e) + 
   \mathcal O( h(x)^{N} ),\qquad h(x)\to0,
\end{equation} 
holds true for any smooth function $f\in C^\infty(\G)$ 
and with $h(x)$ the geodesic distance from $x$ to 
the identity element $e$. An explicit construction of
operators $\partial_x^{(\alpha)}$ in terms of $q^\alpha(x)$
can be found in \cite[Section 10.6]{RT-book}.
In addition to these differential operators
$\partial_x^{(\alpha)}\in \Diff^{|\alpha|}(\G)$ we introduce
operators $\partial_x^{\alpha}$ as follows. Let
$\partial_{x_j}\in \Diff^{1}(\G)$, $j=1,\ldots,n$, be a 
collection of left invariant first order differential
operators corresponding to some linearly independent family
of the left-invariant vector fields on $G$. We denote
$\partial_x^{\alpha}=\partial_{x_1}^{\alpha_1}\cdots
\partial_{x_n}^{\alpha_n}$. We note that in most estimates
we can freely replace operators $\partial_x^{(\alpha)}$
by $\partial_x^{\alpha}$ and the other way around since
they can be clearly expressed in terms of each other.

After fixing the notion of difference operators, we have to specify orders. Each of the bi-invariant subspaces $\mathcal H^\xi$ of Theorem~\ref{thmPW} is an eigenspace of the Laplacian $\Lap$ on $\G$ with the
corresponding eigenvalue $-\lambda_\xi^2$. Based on these eigenvalues we define $\langle\xi\rangle = (1+\lambda_\xi^2)^{1/2}$. In particular we recover the familiar characterisation 
$$
   f\in H^s(\G) \qquad\Longleftrightarrow\qquad \langle\xi\rangle^s \widehat f(\xi) \in \ell^2(\widehat\G).
$$
We are now in a position to formulate the main result of this paper.

\begin{thm}\label{thm:main}
Let $A$ be a linear continuous operator from 
$C^\infty(\G)$ to $\mathcal D'(\G)$. 
Then the following statements are equivalent:
\begin{enumerate}
\item[(A)] 
$A\in\Psi^m(\G)$.
\item[(B)] 
For every left-invariant differential operator 
$P_x\in\Diff^k(\G)$ of order $k$ and every difference operator 
$Q_\xi\in\diff^\ell(\widehat\G)$ of order $\ell$ the symbol estimate
$$
  \| Q_\xi P_x  \sigma_A(x,\xi) \|_{op} \le C_{Q_\xi P_x} \langle\xi\rangle^{m-\ell}
$$
is valid.
\item[(C)] 
For any admissible selection 
$\triangle_1,\ldots,\triangle_n \in\diff^1(\widehat\G)$  we have
$$
  \| \triangle_\xi^\alpha \partial_x^{\beta} 
  \sigma_A(x,\xi)\|_{op} \le C_{\alpha\beta} 
  \langle\xi\rangle^{m-|\alpha|}
$$
for all multi-indices $\alpha,\beta\in\mathbb N_0^n$.
Moreover, $\singsupp R_A(x,\cdot) \subseteq \{e\}$.
\item[(D)] 
For a strongly admissible selection 
$\triangle_1,\ldots,\triangle_n \in\diff^1(\widehat\G)$  we have
$$
  \| \triangle_\xi^\alpha \partial_x^{\beta} 
  \sigma_A(x,\xi)\|_{op} \le C_{\alpha\beta} 
  \langle\xi\rangle^{m-|\alpha|}
$$
for all multi-indices $\alpha,\beta\in\mathbb N_0^n$.
\end{enumerate}
\end{thm}
The set of symbols $\sigma_A$ satisfying either of
conditions $(B)$--$(D)$ will be denoted by
$\mathscr S^{m}(G)$.

Among other things, this theorem removes the assumption of the
conjugation invariance from the list of conditions 
used in \cite{RT-book}.
For $u\in \G$,
denote $A_u f:= A(f\circ\phi)\circ\phi^{-1}$,
where $\phi(x)=xu$. It can be shown that
\begin{eqnarray*}
  R_{A_u}(x,y)  =  R_A(xu^{-1},uyu^{-1})\; \textrm{ and }\;
  \sigma_{A_u}(x,\xi)  =  \xi(u)^\ast\ \sigma_A(xu^{-1},\xi)\ \xi(u).
\end{eqnarray*}
On one hand, we have $A\in\Psi^m(G)$ if and only if 
$A_u\in\Psi^m(G)$ for all $u\in G$. On the other hand,
as a corollary of Theorem \ref{thm:main} this can be expressed
in terms of difference operators, showing that the difference
conditions are conjugation invariant. In fact, to draw
such a conclusion, we even do not need to know that
conditions $(B)$--$(D)$ characterise the class $\Psi^m(G)$:
\begin{cor}\label{COR:RTW-conj-inv}
Let the symbol $\sigma_A$ of a linear continuous operator
from $C^\infty(G)$ to $\mathcal D'(\G)$ satisfy either of
the assumptions $(B)$, $(C)$, $(D)$ of Theorem
\ref{thm:main}. Then the symbols $\sigma_{A_u}$ satisfy the
same symbolic estimates for all $u\in G$.
\end{cor} 
This statement follows immediately from the condition 
$(B)$ if we observe that the difference operators
applied to the symbol $\sigma_{A_u}$ just lead to a different
set of difference operators in the symbolic inequalities.
Alternatively, one can observe that the change of variables
by $\phi$ amounts to the change of the basis in the representation
spaces of $\Gh$, thus leaving the condition $(B)$ invariant
again.

\begin{ex}
On the torus $\Tn=\Rn/\Zn$, the family of
functions $q_j(x)=\erm^{2\pi\irm x_j}-1$, $j=1,\ldots,n$,
gives rise to a strongly admissible collection of
difference operators $$\triangle_j a(\xi)=
\triangle_{q_j} a(\xi)=a(\xi+e_j)-a(\xi),\qquad
j=1,\ldots,n,$$ with $\xi\in \widehat{\Tn}\simeq\Zn$,
where $e_j$ is the $j^{th}$ unit vector in $\Zn$.
As an immediate consequence of Theorem \ref{thm:main}
we recover the fact that the class $\Psi^m(\Tn)$ can be
characterised by the difference conditions on their
toroidal symbols:
$$
  | \triangle_\xi^\alpha \partial_x^{\beta} 
  \sigma_A(x,\xi)| \le C_{\alpha\beta} 
  (1+|\xi|)^{m-|\alpha|}
$$
uniformly for all $x\in\Tn$ and $\xi\in\Zn$.
We refer to \cite{RT-torus} for details of the corresponding
toroidal quantisation of operators on $\Tn$. See also \cite{Agranovich2}.
\end{ex} 

In Theorem~\ref{thm:main} we characterise H\"ormander 
operators entirely by symbol assumptions based on a set 
of admissible difference operators. There does not seem 
to be a canonic choice for difference operators, for 
different applications different selections of them seem 
to be most appropriate. We will comment on some of them.
\begin{enumerate}
\item {\sl Simplicity.} Difference operators have a simple 
structure, if the corresponding functions are just matrix 
coefficients of irreducible representations. 
Then application of difference operators at fixed 
$\xi$ involves only matrix entries from finitely many 
(neighbouring) representations.
\item {\sl Taylor's formula.} Any admissible selection of 
first order differences allows for a Taylor series expansion. 
\item {\sl Leibniz rule.} Leibniz rules have a particularly 
simple structure if one chooses difference operators of
a certain form, see Proposition \ref{prop:Ga-Leibniz}.
Moreover, operators of this form yield a strongly 
admissible collection, see Lemma \ref{LEM:RTW-strong-adms}.
\end{enumerate}

\begin{ex}\label{EX:RTW-SU2}
We will conclude this section with an application
for the particular group $\SU2$. Certain explicit 
calculations on $\SU2$ have been done in \cite{RT-groups} and 
the background on the
necessary representation theory of $\SU2$ can be also found in 
\cite{RT-book}. Representations on $\SU2$ are parametrised by 
half-integers $\ell\in\frac12\mathbb N$ 
(so-called quantum numbers) and are of dimension $
d_\ell = 2\ell+1$. 
Thus, irreducible representations of $\SU2$ are given by matrices
 $t^\ell(x)\in\C^{(2\ell+1)\times(2\ell+1)}$,
$\ell\in\frac12\N$, after some
choice of the basis in the representation spaces.
One admissible selection of difference operators corresponds
to functions $q_+, q_0,q_0\in C^\infty(\SU2)$ defined by
\begin{eqnarray*}
  q_- = t^{1/2}_{-1/2,+1/2}, \quad
  q_+ = t^{1/2}_{+1/2,-1/2}, \quad
  q_0 = 
  t^{1/2}_{-1/2,-1/2}-t^{1/2}_{+1/2,+1/2},
\end{eqnarray*}
where we denote
\begin{eqnarray*}
  t^{1/2}  = 
  \begin{pmatrix}
    t^{1/2}_{-1/2,-1/2} & t^{1/2}_{-1/2,+1/2} \\
    t^{1/2}_{+1/2,-1/2} & t^{1/2}_{+1/2,+1/2}
    \end{pmatrix}.
\end{eqnarray*}
Our analysis on $\SU2$ is in fact equivalent to the 
corresponding analysis on the $3$-sphere $\S3\cong{\rm SU}(2)$,
the isomorphism given by the identification of $\SU2$ with
$\S3\subset {\mathbb H}$ in the quaternion space
${\mathbb H}$, with the quaternionic product
on $\S3$ corresponding to the matrix multiplication on $\SU2$. 
Writing explicitly an isomorphism 
$\Phi:\S3\to\SU2$ as
\begin{equation*}
  (x_0,x_1,x_2,x_3)=x\mapsto
  \Phi(x) = \begin{pmatrix}
        x_0 + {\rm i} x_3 & x_1 + {\rm i} x_2 \\
        -x_1+ {\rm i} x_2 & x_0 - {\rm i} x_3
  \end{pmatrix},
\end{equation*}
the identity matrix $e\in\SU2$ corresponds to the
vector ${\bf 1}$ in the basis decomposition
$
  x=(x_m)_{m=0}^3\mapsto x_0{\bf 1} + x_1{\bf i} + 
  x_2 {\bf j} + x_3 {\bf k}
$
on $\mathbb H$.
We note that the family $\triangle_+:=\triangle_{q_+}$,
$\triangle_-:=\triangle_{q_-}$ and
$\triangle_0:=\triangle_{q_0}$ is not strongly admissible
because in addition to $e\in\SU2$ 
(or to ${\bf 1}\in\mathbb H$) they have another common
zero at $-e$ (or at $-{\bf 1}$). 

Following \cite{RT-book} and \cite{RT-groups}, we simplify
the notation on $\SU2$ and $\S3$ by writing 
$\sigma_A(x,\ell)$ for $\sigma_A(x,t^\ell)$, 
$\ell\in\frac12\N$, and we refer to these works for the
explicit formulae for the difference operators
$\triangle_+, \triangle_-, \triangle_0$.
We denote $\triangle_\ell^\alpha=\triangle_+^{\alpha_1}
\triangle_-^{\alpha_2}\triangle_0^{\alpha_3}$.
In \cite{RT-groups} we proved that the operator and the
Hilbert-Schmidt norms are uniformly equivalent for symbols of
pseudo-differential operators from $\Psi^m(\SU2)$.
As a corollary to
Theorem \ref{thm:main}, Corollary
\ref{COR:RTW-conj-inv}, and \cite[Theorem 12.4.3]{RT-book}
we get
\begin{thm}\label{THM:SU2thm}
Let $A$ be a linear continuous operator from 
$C^\infty(\SU2)$ to $\mathcal D'(\SU2)$. Then
$A\in\Psi^m(\SU2)$ if and only if
$\singsupp R_A(x,\cdot) \subset \{e\}$ and
$$
  \| \triangle_\ell^\alpha \partial_x^{\beta} 
  \sigma_A(x,\ell)\|_{HS} \le C_{\alpha\beta} 
  (1+\ell)^{m-|\alpha|}
$$
for all multi-indices $\alpha,\beta\in\mathbb N^3$
and all $l\in\frac12\N$.
Moreover, in these cases we also have the rapid
off-diagonal decay property of symbols, namely, we have
\begin{equation*}\label{Sigmainequalities}
  \left| \triangle_\ell^\alpha \partial_x^\beta 
  \sigma_{A}(x,\ell)_{ij}\right|
  \leq C_{A\alpha\beta m N}\ (1+ |i-j|)^{-N}
  \ (1+\ell)^{m-|\alpha|}
\end{equation*}
uniformly in $x\in\SU2$, for every $N\geq 0$,
all $\ell\in\frac12\N$,
every multi-indices $\alpha,\beta\in\mathbb N_0^3$,
and for all matrix column/row numbers $i,j$.
\end{thm} 
We note that Theorem \ref{THM:SU2thm} holds in exactly
the same way if
we replace $\SU2$ by $\S3$. Moreover, according
to Theorem \ref{thmPW} the statement of 
Theorem \ref{THM:SU2thm} holds without the kernel
condition $\singsupp R_A(x,\cdot) \subset \{e\}$
provided that instead of the admissible family 
$\triangle_+, \triangle_-, \triangle_0$ we take
a strongly admissible family of difference operators
in $\diff^1(\widehat{\SU2})$.
\end{ex}


\section{Proof of Theorem~\ref{thm:main}}\label{sec3}

\subsection{$(A)\implies (C)$} By \cite[Thm.~10.9.6]{RT-book}
(and in the notation used there) we know that 
condition (A) is equivalent to
$\sigma_A\in\Sigma^m(\G) = \bigcap_k \Sigma_k^m(\G)$, 
where $\Sigma_0^m(\G)$ already corresponds to assumption (C). 
There is nothing to prove.

\subsection{$(C)\implies (D)$} Evident.

Also the following partial converse is true. 
If (D) is satisfied, then the symbol estimates 
imply for the corresponding  right-convolution 
kernel $R_A(x,\cdot)$ of $A$ that
$$
\singsupp R_A(x,\cdot)\subseteq \bigcap_{j=1}^n\{ y\in\G : q_j(y)=0\}=\{e\}.
$$ 

\subsection{$(D)\implies(B)$} 
For a given strongly  admissible selection of 
first order differences, we can apply Taylor's 
formula to the function $q_Q$ corresponding to the 
difference operator $Q\in\diff^\ell(\widehat\G)$. Hence, we obtain

\begin{lem} Let $Q\in\diff^\ell(\widehat\G)$ 
be an arbitrary difference operator of order $\ell$.
Then 
$$
  Q = \sum_{\ell\le |\alpha|\le N-1} c_\alpha 
  \triangle_\xi^\alpha + Q_N
$$
for suitable constants $c_\alpha\in\C$ and a 
difference operator $Q_{N}\in\diff^{N}(\widehat\G)$. 
\end{lem}

Consequently, (D) together with (B) for differences 
of order larger than $N$ implies (B). More precisely, 
with the notation $\|\sigma(x,\xi)\|_{C^\infty}$ for 
`any' $C^\infty$-seminorm of the form
$\sup_x \| P_x \sigma(x,\cdot)\|_{op}$  and
$$
  \mathscr R^m_N = \{ \sigma_A : 
  \;\|Q_\xi \sigma_A(x,\xi)\|_{C^\infty} \lesssim 
  \langle\xi\rangle^{m-\ell}\quad\forall \ell>N\, 
  \ \forall Q_\xi \in\diff^\ell(\widehat\G) \}
$$
and
$$
  \mathscr S^m_\ell = \{ \sigma_A : 
    \;\|\triangle_\xi^\alpha \sigma_A(x,\xi)\|_{C^\infty} 
    \lesssim \langle\xi\rangle^{m-\ell} \quad\forall |\alpha|=\ell\},
$$
we obtain
$$
  \mathscr R^m_0 = \mathscr S^m_0 \cap \mathscr S^m_1 \cap 
  \cdots \cap \mathscr S^m_{N} \cap \mathscr R^m_N, 
$$
while (B) corresponds to taking all symbols from 
$\mathscr R^m_0$. This implies
$$
  \mathscr R^m_0 = \left(\bigcap_{\ell\ge 0} \mathscr S^m_\ell\right)  \cap \mathscr R^m_\infty 
$$
with
$$
   \mathscr R_\infty^m = \{ \sigma_A : 
   \|Q_\xi\sigma(x,\xi)\|_{C^\infty} \lesssim  C_N 
   \langle\xi\rangle^{-N} \quad\forall N\, 
   \forall Q_\xi\in\diff^\infty(\widehat\G)\}.
$$
Assumption (D) corresponds to the first intersection 
$\mathscr S^m = \bigcap_\ell \mathscr S^m_\ell$.
 It remains to show that 
 $\mathscr S^m\subset\mathscr R^m_\infty $, 
 which is done in the sequel.

Let $Q_\xi\in\diff^\infty(\widehat\G)$, 
i.e. $q_Q\in C^\infty(\G)$ is vanishing to infinite 
order at $e$. Let further $p(y) = \sum_{j=1}^n q_j^2(y)$. 
Then $p(y)\ne0$ for $y\ne e$ and $p$ vanishes of the second 
order at $y=e$. Hence,
$q_Q(y) / |p(y)|^N\in C^\infty(\G)$ for arbitrary $N$ and, 
therefore,
$$
q_Q(y) R_A(x,y) = \frac{q_Q(y)}{|p(y)|^N} |p(y)|^N R_A(x,y)
$$
can be estimated in terms of symbol bounds for 
$\sigma_A\in\mathscr S^m$. Indeed,  by Lemma~\ref{lem32}
it follows that 
$|p(y)|^N R_A(x,y) \in H^{-m-c+2N}(\G\times\G)$. 
Letting $N\to\infty$ proves the claim.

\begin{lem}\label{lem32}
Let $\sigma_A\in \mathscr S^{m}$. Then 
$R_A(x,\cdot)\in H^{-m-c}(\G)$ for all $c>\frac12\dim \G$.
\end{lem}
\begin{proof} 
This just follows from 
$\|\sigma_A(x,\xi)\|_{HS} \le \sqrt{d_\xi} 
\|\sigma_A(x,\xi)\|_{op}$ in combination with 
the observation that 
$$ \sum_{[\xi]} d_\xi^2 \langle\xi\rangle^{-2c}< \infty $$
for $c > \frac12\dim\G$ (which is equivalent to 
$H^c(\G)\hookrightarrow L^2(\G)$ being Hilbert-Schmidt). 
\end{proof}

\subsection{$(B)\Longleftrightarrow(C)\Longleftrightarrow(D)$} 
Evident from the above. 

\subsection{$(B) \implies (A)$} We use the commutator 
characterisation of the class
$\Psi^m(\G)$, \cite[Thm.~10.7.7]{RT-book}. 
It says that it is sufficient to show that
\begin{enumerate}
\item any symbol $\sigma_A$ satisfying condition (B) 
gives rise to a bounded linear operator 
$A: H^m(\G)\to L^2(\G)$; 
\item the symbol $\sigma_{[X,A]}$ of the commutator of 
$A$ and any left-invariant vector field $X$ also
satisfies (B).
\end{enumerate}
The nuclearity of $C^\infty(\G)$ allows one 
to split the first statement into the corresponding results 
for purely $x$- and purely $\xi$-dependent symbols;  
both situations are evident.  
As for the second statement, we use the symbolic 
calculus established in \cite[Thm.~10.7.9]{RT-book} 
for classes of symbols including, in particular, the class
$\mathscr S^m$, which writes the symbol of the composition as
the asymptotic sum
$$
 \sigma_{[X,A]} \sim \sigma_X \sigma_A - \sigma_A\sigma_X + 
 \sum_{|\alpha|>1} \frac1{\alpha!} (\triangle_\xi^\alpha \sigma_X)
(\partial_x^{(\alpha)}\sigma_A).
$$
The first term $\sigma_X \sigma_A - \sigma_A\sigma_X$
is on the level of right-convolution kernels 
$R_A(x,\cdot)$ given by
$$
   (X_y  - X^R_y) R_A(x,y),
$$ 
where $X^R$ is the right-invariant vector field tangent to $X$ 
at $e$ (and the index $y$ describes the action as a 
differential operator with respect to the $y$-variable). 
Apparently, $X_y-X^R_y$ is a smooth vector field on $\G$ and
can thus be written as 
$$
   \sum_{j=1}^n q_j(y) \partial_j^{R}
$$ 
with $q_j(e)=0$ and $\partial_j^R$ suitable 
right-invariant derivatives. Hence, on the symbolic 
level $\sigma_X \sigma_A - \sigma_A\sigma_X = \sum_{j=1}^n 
\triangle_j (\sigma_A\sigma_{\partial_j})$ with 
$\triangle_j\in\diff^1(\widehat\G)$, defined by
$\triangle_j\widehat{f}(\xi)=\widehat{q_j f}(\xi)$.
Now the statement follows from 
Lemma \ref{LEM:DOps}
by the aid of the asymptotic 
Leibniz rule \cite[Thm.~10.7.12]{RT-book}:
\begin{lem}\label{LEM:RTW-ass-Leibniz}
Let  $\triangle_1,\ldots,\triangle_n\in\diff^1(\widehat\G)$ 
be a set of admissible differences. 
For any $Q_\xi\in\diff^\ell(\widehat \G)$ 
there exist $Q_{\xi,\alpha}\in\diff^\ell(\widehat\G)$ such that
  $$
     Q_\xi( \sigma_1\sigma_2) \sim (Q_\xi \sigma_1)\sigma_2 
     +\sum_{|\alpha|>0} (Q_{\xi,\alpha} \sigma_1) 
     (\triangle_\xi^\alpha \sigma_2)
  $$
  for any (fixed) set of admissible differences.
\end{lem}

\begin{lem}\label{LEM:DOps}
Let $A$ be a differential operator of order $m$. 
Then $\sigma_A\in \mathscr S^m$.
\end{lem}

\begin{proof}
Lemma \ref{LEM:DOps} is a consequence of the already
proved implications
$(A)\Longrightarrow(C)\Longrightarrow(D)\Longrightarrow(B)$.
Explicit formulae for the difference operators applied to 
symbols of differential operators can be also found
in \cite[Prop.~10.7.4]{RT-book} and provide an elementary proof of this statement.
\end{proof}


\section{Ellipticity and Leibniz formula}\label{sec4}

As an application of Theorem~\ref{thm:main} 
we will give a characterisation of the elliptic 
operators in $\Psi^m(\G)$ in terms of their global symbols.

\begin{thm}\label{THM:RTW-ellipticity}
An operator $A\in\Psi^m(\G)$ is elliptic if and only if
its matrix valued symbol $\sigma_A(x,\xi)$ is invertible 
for all but finitely many $[\xi]\in\Gh$, and for all such $\xi$
satisfies
\begin{equation}\label{EQ:RTW-ellsymbs1}
  \|\sigma_A(x,\xi)^{-1}\|_{op} \le C  \langle\xi\rangle^{-m}
\end{equation} 
for all $x\in\G$.
\end{thm}
Thus, both statements are equivalent to the existence of $B\in\Psi^{-m}(G)$
such that $I-BA$ and $I-AB$ are smoothing smoothing.
For the ellipticity condition \eqref{EQ:RTW-ellsymbs1} on the
general matrix level it is not enough to assume that
$|\det\sigma_A(x,\xi)|^{1/d_\xi}\ge C \langle\xi\rangle^{m}$ 
due to the in general growing dimension of the matrices. 
However, if we assume that the smallest singular value of 
the matrix $\sigma_A(x,\xi)$ is 
greater or equal than $C\jp{\xi}^{m}$ uniformly in $x$ 
and (all but finitely many) $\xi$,
then condition \eqref{EQ:RTW-ellsymbs1} follows.

Let the collection $q_1,\ldots,q_n$ give an admissible
collection of difference operators and let
$\partial_x^{(\gamma)}$ be the corresponding family of 
differential operators as in the Taylor
expansion formula \eqref{EQ:RTW-Taylor-exp}.
As an immediate corollary of Theorem
\ref{THM:RTW-ellipticity} and \cite[Thm.~10.9.10]{RT-book}
we get
\begin{cor}
Let $A\in\Psi^m(G)$ be elliptic.
Moreover, assume that 
$
  \sigma_A(x,\xi) \sim \sum_{j=0}^\infty \sigma_{A_j}(x,\xi),
$
where $\Op(\sigma_{A_j})\in \Psi^{m-j}(G)$. Let
$
  \sigma_B(x,\xi) \sim \sum_{k=0}^\infty \sigma_{B_k}(x,\xi),
$
where $\sigma_{B_0}(x,\xi)=\sigma_{A_0}(x,\xi)^{-1}$ for 
large $\xi$, and the
symbols $\sigma_{B_k}$ 
are defined recursively by
\begin{equation}\label{EQ:Gparametrix1}
  \sigma_{B_N}(x,\xi)
  =  -\sigma_{B_0}(x,\xi)
  \sum_{k=0}^{N-1} \sum_{j=0}^{N-k} \sum_{|\gamma|=N-j-k}
  \frac{1}{\gamma!} 
  \left[ \triangle_\xi^\gamma \sigma_{B_k}(x,\xi) \right]
  \partial_x^{(\gamma)} \sigma_{A_j}(x,\xi).
\end{equation}
Then $\Op(\sigma_{B_k})\in \Psi^{-m-k}(G)$, $B=\Op(\sigma_B)\in\Psi^{-m}(G)$,
and the operators  $AB-I$ and $BA-I$ are in $\Psi^{-\infty}(G)$.
\end{cor} 

First, we give some preliminary results.
In Lemma \ref{LEM:RTW-ass-Leibniz} we gave
an asymptotic Leibniz formula 
but here we will present its finite version.
Given a continuous unitary matrix representation
$\xi=\begin{bmatrix} \xi_{ij} \end{bmatrix}_{1\leq i,j\leq \ell}:
G\to\mathbb C^{\ell\times \ell}$,  
$\ell=d_\xi$,
let $q(x)=\xi(x)-I$ (i.e. $q_{ij}=\xi_{ij}-\delta_{ij}$
with Kronecker's deltas $\delta_{ij}$), and define
$$
  \D_{ij}\widehat{f}(\xi) := \widehat{q_{ij} f}(\xi).
$$
In the previous notation, we can write $\D_{ij}=\Delta_{q_{ij}}$.
However, here, to emphasize that $q_{ij}$'s correspond to the
same representation, and to distinguish with difference
operators with a single subindex, we will use the notation
$\D$ instead.
For a multi-index $\gamma\in \N_0^{\ell^2}$, we will write
$|\gamma|=\sum_{i,j=1}^\ell|\gamma_{ij}|$, and for higher order
difference operators we write
$\D^\gamma=\D_{11}^{\gamma_{11}}\D_{12}^{\gamma_{12}}\cdots
\D_{\ell,\ell-1}^{\gamma_{\ell,\ell-1}}\D_{\ell\ell}^{\gamma_{\ell\ell}}$.

For simplicity, we will abbreviate writing $a=a(\xi)$ and $b=b(\xi)$.
Let $k\in\N_0$ and let $\alpha,\beta\in\N^k$ be
such that $1\leq \alpha_j,\beta_j\leq \ell$ for all $j\in\{1,\cdots,k\}$.
Let us define a {\it grand $k$th order difference} 
$\Darr^k$ by
$$
  (\Darr^k a)_{\alpha\beta} = \Darr^k_{\alpha\beta} a :=
  \D_{\alpha_1\beta_1} \cdots \D_{\alpha_k\beta_k} a,
$$
where we may note that the first order differences 
$\D_{\alpha_j\beta_j}$ commute with each other.
We may compute with ``matrices''
$$
  \Darr^k a =
  \begin{bmatrix} \Darr^k_{\alpha\beta} a
  \end{bmatrix}_{\alpha,\beta\in\{1,\cdots,\ell\}^k}
$$
using the natural operations, e.g.
$$
  \left((\Darr^k a)(\Darr^k b)\right)_{\alpha\beta}
  := \sum_{\gamma\in\{1,\cdots,\ell\}^k} 
  \left(\Darr^k_{\alpha\gamma} a\right)
  \left(\Darr^k_{\gamma\beta} b\right).
$$
In contrast to Lemma \ref{LEM:RTW-ass-Leibniz} with an
asymptotic Leibniz formula for arbitrary
difference operators,
operators $\Darr^k$ satisfy the finite Leibniz formula:
\begin{prop}\label{prop:Ga-Leibniz}
For all $k\in\N$ and $\alpha,\beta\in\{1,\cdots,\ell\}^k$ we have
$$
\Darr^k_{\alpha\beta}(ab)=
\sum_{|\varepsilon|,|\delta|\leq k\leq |\varepsilon|+|\delta|}
C^k_{\alpha\beta\varepsilon\delta}\ (\D^\varepsilon a)\ 
(\D^\delta b),
$$
with the summation taken over all 
$\varepsilon,\delta\in\N_0^{\ell^2}$ satisfying
$|\varepsilon|,|\delta|\leq k\leq |\varepsilon|+|\delta|$. 
In particular, for
$k=1$, we have
\begin{eqnarray}\label{EQ:G-leinbiz-1}
\D_{ij}(ab) & = & \left(\D_{ij} a\right) b + a\left(\D_{ij} b\right)
  + \left( (\Darr^1 a)(\Darr^1 b) \right)_{ij} \\ \nonumber 
  & = & 
  \left(\D_{ij} a\right) b
  + a \left(\D_{ij} b\right)
  + \sum_{k=1}^n
  \left( \D_{ik} a \right)
  \left( \D_{kj} b \right).
\end{eqnarray} 
\end{prop}
\begin{proof}
Let us first show the case of $k=1$ given
in formula \eqref{EQ:G-leinbiz-1}.
Recalling the notation $q_{ij}=\xi_{ij}-\delta_{ij}$, 
from the unitarity of the
representation $\xi$, for all $x,y\in G$, we readily obtain the 
identity
\begin{equation}\label{EQ:G-leibnitz-2}
q_{ij}(x)=q_{ij}(x y^{-1}) + q_{ij}(y)+\sum_{k=1}^n
q_{ik}(xy^{-1})q_{kj}(y).
\end{equation} 
If now $f,g\in C^\infty(G)$, we consequently get
$$
q_{ij}\cdot(f*g)=(q_{ij}f)*g+f*(q_{ij}g)+\sum_{k=1}^n
(q_{ik}f)*(q_{kj}g),
$$
which on the Fourier transform side gives \eqref{EQ:G-leinbiz-1}.

For $k\geq 2$, let us argue inductively.
Taking differences is linear, so the main point is to analyse
the application of a first order difference
$\D_{\alpha_1\beta_1}$ to the term
$\left( (\Darr^1 a)(\Darr^1 b) \right)_{\alpha_2\beta_2}$.
By \eqref{EQ:G-leinbiz-1} we have
\begin{eqnarray*}
  \D_{\alpha_1\beta_1}
  \left( (\Darr^1 a)(\Darr^1 b) \right)_{\alpha_2\beta_2}
  & = &
  \left( (\Darr^1 \D_{\alpha_1\beta_1} a) (\Darr^1 b)
  \right)_{\alpha_2\beta_2}
  + \left( (\Darr^1 a) (\Darr^1\D_{\alpha_1\beta_1} b)
  \right)_{\alpha_2\beta_2} \\
  && + \left( (\Darr^2 a)(\Darr^2 b) \right)_{\alpha\beta},
\end{eqnarray*}
where $\alpha=(\alpha_1,\alpha_2)$ and $\beta=(\beta_1,\beta_2)$.
Proceeding by induction, we notice that
applying a first order difference to the term
$\left( (\Darr^k a)(\Darr^k b) \right)_{\alpha\beta}$
would introduce a sum of old type terms plus the new term
$\left( (\Darr^{k+1} a)(\Darr^{k+1} b)
\right)_{\widetilde\alpha\widetilde\beta}$.
Thus, a $k$th order difference operator
applied to the product symbol $ab$
produces a linear combination of terms of the form
$
  (\D^\varepsilon a)(\D^\delta b),
$
where
$
  |\varepsilon|,|\delta|\leq k\leq |\varepsilon|+|\delta|,
$
completing the proof.
\end{proof}

Difference operators $\D$ associated to the representations
may be particularly useful in view of the finite Leibniz
formula in Proposition \ref{prop:Ga-Leibniz} and the fact that
the collection of all differences $\D_{ij}$ over all $\xi$ and
$i,j$ is strongly admissible:
\begin{lem}\label{LEM:RTW-strong-adms}
The family of difference operators associated to the family of
functions
$\{q_{ij}=\xi_{ij}-\delta_{ij}\}_{[\xi]\in\Gh,\ 1\leq i,j\leq d_\xi}$
is strongly admissible. Moreover, this family has a finite
subfamily 
associated to finitely many representations which is still 
strongly admissible.
\end{lem} 
\begin{proof}
We observe that there exists a homomorphic 
embedding of $G$ into $U(N)$ 
for a large enough $N$ and this embedding itself 
is a representation of
$G$ of dimension $N$.  Decomposing this representation into 
irreducible components gives a finite collection of representations.
The common zero set of the corresponding family $\{q_{ij}\}$ is $e$
which means that it is strongly admissible.
\end{proof} 
We note that on the group $\SU2$ or on $\S3$, 
by this argument, or by the 
discussion in Section \ref{sec2} the four function $q_{ij}$
corresponding to the representation $[t^\ell]\in\widehat{\SU2}$ 
with $\ell=\frac12$ 
of dimension two, $d_\ell=2$, already give a strongly admissible
collection of four difference operators.
We can now apply Proposition \ref{prop:Ga-Leibniz}
to the question of inverting the symbols on $G\times\Gh$.
We formulate this for symbol classes with $(\rho,\delta)$
behaviour as this will be also used in Section
\ref{sec5}.
  
\begin{lem}\label{LEM:RTW-inv-inverse}
Let $m\geq m_0$, $1\ge \rho>\delta\ge 0$ and let 
us fix difference operators 
$\{\D_{ij}\}_{1\leq i,j\leq d_{\xi_0}}$
corresponding to some representation $\xi_0\in\Gh$.
Let the matrix symbol $a=a(x,\xi)$ 
satisfy
\begin{equation}\label{EQ:RTW-a-ineqs}
 \| \D^\gamma  \partial_x^\beta
  a(x,\xi)\|_{op} \le C_{\beta\gamma} 
  \langle\xi\rangle^{m-\rho|\gamma|+\delta |\beta|}
\end{equation} 
for all multi-indices $\beta,\gamma$.
Assume also that  
$a(x,\xi)$ is invertible for all $x\in G$ and $[\xi]\in\Gh$,
and satisfies
\begin{equation}\label{EQ:G-ell-1}
\|a(x,\xi)^{-1}\|_{op}\leq C\jp{\xi}^{-m_0}
\end{equation} 
for all $x\in G$ and $[\xi]\in\Gh$, 
and if $m_0\ne m$ in addition that
\begin{equation}\label{EQ:hypo}
   \n{a(x,\xi)^{-1} 
   \br{\D^\gamma \partial_x^\beta a(x,\xi)}
   }_{op} \le C  \langle\xi\rangle^{-\rho |\gamma|+\delta|\beta|}
\end{equation} 
for all $x\in G$ and $[\xi]\in\Gh$. Then the matrix symbol
$a^{-1}$ defined by $a^{-1}(x,\xi)=a(x,\xi)^{-1}$ satisfies
\begin{equation}\label{EQ:RTW-a-ineqs2}
 \| \D^\gamma  \partial_x^\beta
  a^{-1}(x,\xi)\|_{op} \le C_{\beta\gamma} 
  \langle\xi\rangle^{-m_0-\rho |\gamma|+\delta |\beta|}
\end{equation} 
for all multi-indices $\beta,\gamma$.
\end{lem} 
\begin{proof}
Let us denote $b(x,\xi):=a(x,\xi)^{-1}$
and estimate $\partial_x^\beta b$ first. 
Suppose we have proved that
\begin{equation}\label{INEQ:inductionhypothesis2}
  \|\partial_x^\beta b(x,\xi)\|_{op}
  \leq C_\beta \langle\xi\rangle^{-m_0+\delta|\beta|}
\end{equation}
whenever $|\beta|\leq k$.
We proceed by induction. Let us study the order $k+1$ 
cases $\partial_{x_j}\partial_x^\beta b$ with
$|\beta|=k$. Since $a(x,\xi)b(x,\xi)=I$, by the usual
Leibniz formula we get
$$
 a\ \partial_{x_j}\partial_x^\beta b =
 -\sum_{\beta_1+\beta_2=\beta+e_j,|\beta_2|\leq|\beta| }
 C_{\beta_1\beta_2}\ (\partial_x^{\beta_1} a)\ 
 (\partial_x^{\beta_2} b) .
$$
From \eqref{EQ:RTW-a-ineqs}, \eqref{EQ:G-ell-1}, \eqref{EQ:hypo} and
\eqref{INEQ:inductionhypothesis2} we obtain the desired 
estimate
$$\|\partial_{x_j}\partial_x^\beta b(x,\xi)\|_{op}
  \leq C_\beta \langle\xi\rangle^{-m_0+\delta|\beta|+\delta}.$$

Let us now estimate $\D^\gamma b$. 
The argument is more complicated than that for
the $\partial_x$-derivatives because the Leibniz formula
in Proposition \ref{prop:Ga-Leibniz} has more terms. 
Indeed, the Leibniz formula in Proposition \ref{prop:Ga-Leibniz}
applied to $a(x,\xi) b(x,\xi)=I$ gives
$$
\D_{ij}b+b(\D_{ij}a) b+
\sum_{k=1}^{d_{\xi_0}} b(\D_{ik}a)\D_{kj}b=0.
$$
Writing these equations for all $1\leq i,j\leq d_{\xi_0}$
gives a linear system on $\{\D_{ij}b\}_{ij}$
with coefficients of the form $I+b\D_{kl}a$
for suitable sets of indices $k,l$. Since
$\|b\D_{kl} a\|_{op}\leq C\jp{\xi}^{-\rho}$, we
can solve it for $\{\D_{ij}b\}_{ij}$
for large $\jp{\xi}$. The inverse of this system
is bounded, and $\|b(\D_{ij}a) b\|_{op}\leq C\jp{\xi}^{-m_0-\rho}$,
implying that $\|\D_{ij} b\|_{op}\leq C\jp{\xi}^{-m_0-\rho}$.

Suppose now we have proved that
\begin{equation}\label{INEQ:inductionhypothesis}
  \|\D^\gamma b(x,\xi)\|_{op}
  \leq C_\gamma \langle\xi\rangle^{-m_0-\rho|\gamma|}
\end{equation}
whenever $|\gamma|\leq k$.
We proceed by induction.
Let us study the order $k+1$ cases $\D_{ij}\D^\gamma b$ with
$|\gamma|=k$.
Now $a(x,\xi)b(x,\xi)=I$, so that
\begin{equation*}\label{EQ:zerotriangleAB}
\D_{ij}\D^\gamma(ab)=0,
\end{equation*}
where by Proposition \ref{prop:Ga-Leibniz}
the right hand side is a sum of terms of the form
$ 
  \left(\D^\varepsilon a\right)\left(\D^\eta b\right),
$
where $|\varepsilon|,|\eta|\leq k+1\leq |\varepsilon|+|\eta|$.
Especially, we look at the terms
$  
\D^\varepsilon a\ \D^\eta b
$
with $|\eta|=k+1$, since
the other terms can be estimated by $\langle\xi\rangle^{-\rho(k+1)}$
due to \eqref{INEQ:inductionhypothesis} and
\eqref{EQ:RTW-a-ineqs}. 
Writing the linear system of equations on $\D_{ij}\D^{\eta'} b$
produced by the Leibniz formula
for all $i,j$ and all $|\eta'|=k$,
we see that the matrix coefficients in front of matrices
$\D_{ij}\D^{\eta'} b$ are sums of the terms of the form
$\D^\varepsilon a$. The main term in each coefficient
is $a(x,\xi)$ corresponding to $\varepsilon=0$, while the
other terms corresponding to $\varepsilon\not=0$ can be
estimated by $\jp{\xi}^{m-\rho}$ in view of
\eqref{INEQ:inductionhypothesis} and
\eqref{EQ:RTW-a-ineqs}. It follows that we can solve
$\D_{ij}\D^{\eta'} b$ from this system and the solution
matrix can be estimated by $\jp{\xi}^{-m_0}$ in view of
\eqref{EQ:G-ell-1} since its main term is $a(x,\xi)$.
Therefore, all terms of the type $\D_{ij}\D^{\eta'}$ can be
estimated by
$$
  \left\|\D_{ij}\D^{\eta'} b(x,\xi)\right\|_{op}
    \leq C_{\eta',i,j}\ \langle\xi\rangle^{-m_0-\rho(k+1)},
$$
which is what was required to prove. By combining these
two arguments we obtain estimate  \eqref{EQ:RTW-a-ineqs2}. 
\end{proof}

\begin{proof}[Proof of Theorem \ref{THM:RTW-ellipticity}]
Assume first that $A\in\Psi^m(G)$ is elliptic. 
Then it has a parametrix $B\in\Psi^{-m}(G)$ such that
$AB-I, BA-I\in\Psi^{-\infty}(G)$. By the composition
formulae for the matrix valued symbols in
\cite[Thm.~10.7.9]{RT-book} and Theorem
\ref{thm:main} it follows that
$\Op(\sigma_A\sigma_B)-I\in\Psi^{-1}(G)$. Consequently,
the product $\sigma_A(x,\xi)\sigma_B(x,\xi)$ is an invertible
matrix for all sufficiently large $\jp{\xi}$ and
$\|(\sigma_A\sigma_B)^{-1}\|_{op}\leq C$. 
It follows then also that $\sigma_A(x,\xi)$ is an invertible
matrix for all sufficiently large $\jp{\xi}$.
From this and
the equality $\sigma_A^{-1}(\sigma_A\sigma_B)=\sigma_B$ we
obtain
$$
\|\sigma_A^{-1}\|_{op}=\|\sigma_B (\sigma_A\sigma_B)^{-1}\|_{op}
\leq \|\sigma_B\|_{op}\|(\sigma_A\sigma_B)^{-1}\|_{op}
\leq C\jp{\xi}^{-m}.
$$

Assume now that $A\in\Psi^m(G)$ has the invertible
matrix symbol $\sigma_A$ satisfying 
$\|\sigma_A(x,\xi)^{-1}\|_{op}\leq C\jp{\xi}^{-m}.$
We can disregard representations $\xi$ with bounded $\jp{\xi}$
because they correspond to a smoothing operator.
From Theorem \ref{thm:main} applied with the
difference operators from 
Lemma \ref{LEM:RTW-strong-adms}, we get that
$\Op(\sigma_A^{-1})\in\Psi^{-m}(G)$ because we have
$\|\D^\gamma\partial_x^\beta (\sigma_A(x,\xi)^{-1})\|_{op}
\leq C\jp{\xi}^{-m-|\gamma|}$ by
Lemma \ref{LEM:RTW-inv-inverse} with $m_0=m$, $\rho=1$ and $\delta=0$.
Since $\Op(\sigma_A^{-1})\in\Psi^{-m}(G)$, Theorem 10.9.10
in \cite{RT-book} implies that $A$ has a parametrix
(given by formula \eqref{EQ:Gparametrix1}). This implies
that $A$ is elliptic.
The proof is complete.
\end{proof}  

\section{Global hypoellipticity}
\label{sec5}
We now turn to the analysis of hypoelliptic operators.
We recall that an operator $A$ is globally hypoelliptic on $G$
if $Au=f$ and $f\in C^\infty(G)$ imply $u\in C^\infty(G)$.
We use the notation 
$\mathscr S^m_{\rho,\delta}(\G)$ for the class of
symbols for which the corresponding 
$(\rho,\delta)$-versions
of symbol estimates are satisfied, namely, 
$\sigma_A\in\Smrd$ if 
for a strongly admissible selection 
$\triangle_1,\ldots,\triangle_n \in\diff^1(\widehat\G)$  we have
$$
  \| \triangle_\xi^\alpha \partial_x^{\beta} 
  \sigma_A(x,\xi)\|_{op} \le C_{\alpha\beta} 
  \langle\xi\rangle^{m-\rho|\alpha|+\delta|\beta|}
$$
for all multi-indices $\alpha,\beta\in\mathbb N_0^n$.
Theorem~\ref{thm:main}  
$(B)\Leftrightarrow (C)\Leftrightarrow(D)$ is valid for all  
$1\ge \rho>\delta\ge 0$. 
By $\Op(\mathscr S^m_{\rho,\delta}(\G))$ we denote the
class of all operators $A$ of the form \eqref{EQ:RTW-quant}
with symbols $\sigma_A\in\Smrd$. 
Such operators can be readily seen to be continuous on
$C^\infty(G)$.
In the previous notation
we have $\mathscr S^m(G)=\mathscr S^m_{1,0}(G)$, and
$\Psi^m(G)=\Op(\mathscr S^m_{1,0}(G))$ by 
Theorem~\ref{thm:main}.

\subsection{}
The knowledge of the full
symbols allows us to establish an analogue of the
well-known hypoellipticity result of H\"ormander \cite{Ho}
on $\Rn$, requiring conditions on lower order terms of
the symbol.
The following theorem is a matrix-valued symbol criterion
for (local) hypoellipticity.

\begin{thm}\label{THM:RTW-hypoellipticity}
Let $m\geq m_0$ and $1\ge\rho>\delta\ge 0$.
Let $A\in\Op(\Smrd)$ be a pseudo-differential operator with 
the matrix-valued symbol 
$\sigma_A(x,\xi)\in\Smrd$ which is invertible 
for all but finitely many $[\xi]\in\Gh$, and for all such $\xi$
satisfies
\begin{equation}\label{EQ:RTW-he1}
   \|\sigma_A(x,\xi)^{-1}\|_{op} \le C  \langle\xi\rangle^{-m_0}
\end{equation} 
for all $x\in\G$. Assume also that for a strongly admissible
collection of difference operators we have
\begin{equation}\label{EQ:RTW-he2}
   \n{\sigma_A(x,\xi)^{-1} 
   \br{\triangle_\xi^\alpha\partial_x^\beta \sigma_A(x,\xi)}
   }_{op} \le C  \langle\xi\rangle^{-\rho |\alpha|+\delta |\beta|}
\end{equation} 
for all multi-indices $\alpha,\beta$, all $x\in G$, and all but
finitely many $[\xi]$.
Then there exists an operator 
$B\in\Op(\mathscr S^{-m_0}_{\rho,\delta}(G))$ such
that $AB-I$ and $BA-I$ belong to $\Psi^{-\infty}(G)$.
Consequently, we have 
$$
\singsupp Au = \singsupp u\; \textrm{ for all }\; u\in\mathcal D'(G).
$$
\end{thm} 
\begin{proof}
First of all, we can assume at any stage of the proof
that \eqref{EQ:RTW-he1} and
\eqref{EQ:RTW-he2} hold for all $[\xi]\in\Gh$ since
we can always modify symbols for small $\jp{\xi}$
which amounts to adding a smoothing operator.
Then we observe that if we apply Theorem
\ref{thm:main} with difference operators $\D$ from
Lemma \ref{LEM:RTW-strong-adms}, it follows from
Lemma \ref{LEM:RTW-inv-inverse} that 
$\sigma_A^{-1}\in \mathscr S_{\rho,\delta}^{-m_0}(G)$.
Let us show next that 
$\sigma_A^{-1} 
\triangle_\xi^\alpha\partial_x^\beta \sigma_A\in
\mathscr S_{\rho,\delta}^{-\rho|\alpha|+\delta|\beta|}(G)$.
Differentiating the equality 
\begin{equation}\label{EQ:RTW-aux-1}
\sigma_A (\sigma_A^{-1} 
\triangle_\xi^\alpha\partial_x^\beta \sigma_A)=
\triangle_\xi^\alpha\partial_x^\beta \sigma_A,
\end{equation} 
we get
$$
\partial_{x_j}\p{\sigma_A^{-1} 
\triangle_\xi^\alpha\partial_x^\beta \sigma_A}=
\sigma_A^{-1}
\partial_{x_j}\triangle_\xi^\alpha\partial_x^\beta \sigma_A
-(\sigma_A^{-1}\partial_{x_j}\sigma_A) (\sigma_A^{-1} 
\triangle_\xi^\alpha\partial_x^\beta \sigma_A).
$$
From this and \eqref{EQ:RTW-he2} it follows that
$\|\partial_{x_j}\p{\sigma_A^{-1} 
\triangle_\xi^\alpha\partial_x^\beta \sigma_A}\|_{op}\leq
C\jp{\xi}^{-\rho|\alpha|+\delta|\beta|+\delta}.$ Continuing this argument, we get
that $\|\partial_x^\gamma\p{\sigma_A^{-1} 
\triangle_\xi^\alpha\partial_x^\beta \sigma_A}\|_{op}\leq
C\jp{\xi}^{-\rho|\alpha|+\delta|\beta+\gamma|}$ for all $\gamma$. For differences,
the Leibniz formula in Proposition \ref{prop:Ga-Leibniz}
applied to \eqref{EQ:RTW-aux-1} gives
\begin{multline*}
\D_{ij}(\sigma_A^{-1}\triangle_\xi^\alpha\partial_x^\beta \sigma_A)
+(\sigma_A^{-1}\D_{ij}\sigma_A) (\sigma_A^{-1}
\triangle_\xi^\alpha\partial_x^\beta \sigma_A) \\ +
\sum_k (\sigma_A^{-1}\D_{ik}\sigma_A)\D_{kj}(\sigma_A^{-1}
\triangle_\xi^\alpha\partial_x^\beta \sigma_A) =
\sigma_A^{-1}\D_{ij}\triangle_\xi^\alpha\partial_x^\beta \sigma_A.
\end{multline*}
Writing these equations for all $i,j$
gives a linear system on 
$\{\D_{ij}(\sigma_A^{-1}
\triangle_\xi^\alpha\partial_x^\beta \sigma_A)\}_{ij}$
with coefficients of the form $I+\sigma_A^{-1}\D_{kl}\sigma_A$
for suitable sets of indices $k,l$. Since
$\|\sigma_A^{-1}\D_{kl}\sigma_A\|_{op}\leq C\jp{\xi}^{-\rho}$
by \eqref{EQ:RTW-he2} and Theorem \ref{thm:main}, we
can solve it for matrices
$\{\D_{ij}(\sigma_A^{-1}
\triangle_\xi^\alpha\partial_x^\beta \sigma_A)\}_{ij}$
for large $\jp{\xi}$. Finally, since
$$
\|(\sigma_A^{-1}\D_{ij}\sigma_A) (\sigma_A^{-1}
\triangle_\xi^\alpha\partial_x^\beta \sigma_A)\|_{op}
\leq C\jp{\xi}^{-\rho|\alpha|+\delta|\beta|-\rho},$$
$$ \|\sigma_A^{-1}\D_{ij}
\triangle_\xi^\alpha\partial_x^\beta \sigma_A\|_{op}
\leq C\jp{\xi}^{-\rho|\alpha|+\delta|\beta|-\rho},
$$
we get that
$\|\D_{ij}\p{\sigma_A^{-1} 
\triangle_\xi^\alpha\partial_x^\beta \sigma_A}\|_{op}\leq
C\jp{\xi}^{-\rho|\alpha|+\delta|\beta|-\rho}$ for all $i,j$. An induction
argument similar to the one in the proof of 
Lemma \ref{LEM:RTW-inv-inverse} shows that
$\sigma_A^{-1} 
\triangle_\xi^\alpha\partial_x^\beta \sigma_A\in
\mathscr S_{\rho,\delta}^{-\rho|\alpha|+\delta|\beta|}(G)$.
Let us now denote $\sigma_{B_0}(x,\xi)=
\sigma_{A_0}(x,\xi)^{-1},$
and let
$\sigma_B\sim\sum_{N=0}^\infty \sigma_{B_N}$, where
\begin{equation}\label{EQ:Gparametrix2}
  \sigma_{B_N}(x,\xi) = -\sigma_{B_0}(x,\xi)
  \sum_{k=0}^{N-1} 
  \sum_{|\gamma|=N-k}
        \frac{1}{\gamma!} \left[
          \triangle_\xi^\gamma \sigma_{A}(x,\xi) \right]
        \partial_x^{(\gamma)} \sigma_{B_k}(x,\xi).
\end{equation}
One can readily check that $B$ is a parametrix of $A$
(this is a special case of \cite[Exercise 10.9.12]{RT-book}). We claim
that $\sigma_{B_N}\in \mathscr S_{\rho,\delta}^{-m_0-N}(G).$
We already know
that $\sigma_{B_0}\in \mathscr S_{\rho,\delta}^{-m_0}(G).$ 
Inductively, in view of
$\sigma_A^{-1} 
\triangle_\xi^\alpha\partial_x^\beta \sigma_A\in
\mathscr S_{\rho,\delta}^{-\rho|\alpha|+\delta|\beta|}(G)$,
formula \eqref{EQ:Gparametrix2} implies that
the order of $B_N$ is $-(\rho-\delta)|\gamma|+\textrm{order}(B_k)
=-(\rho-\delta)|\gamma|-m_0-(\rho-\delta)k=-m_0-(\rho-\delta)N.$ Thus, 
$B\in\Op(\mathscr S^{-m_0}_{\rho,\delta}(G))$ is
a parametrix for $A$.

Finally, the pseudo-locality of the operator $A$ 
(which can be proved directly from the symbolic calculus)
implies that $\singsupp Au\subset\singsupp u$. Writing
$u=B(Au)-Ru$ with $R\in\Psi^{-\infty}(G)$, we get that
$$
\singsupp u\subset \singsupp (Au) \bigcup
\singsupp (Ru),
$$
so that $\singsupp Au=\singsupp u$ because
$Ru\in C^\infty(G)$.
\end{proof} 

To obtain global hypoellipiticity, it is sufficient to show that an operator $A\in\Psi(\G)$ has a parametrix $B$ satisfying subelliptic estimates $\|B f\|_{H^s} \lesssim \|f\|_{H^{s+m}}$ for some constant $m$ independent of $s\in\R$. 

\subsection{Examples}
We will conclude with a collection of examples on the group $\SU2$.
These examples can be also viewed as examples on the
3-sphere $\S3$, see Example \ref{EX:RTW-SU2}.
Thus, in what follows we denote $G=\S3\cong{\rm SU}(2)$.
In contrast to \cite{RT-book} we will base considerations here
on the four difference operators $\D_{ij}$, $i,j\in\{1,2\}$, 
defined in terms of the first non-trivial representation 
corresponding to $\ell=\frac12$. 
They are related to $\triangle_\eta$, $\eta\in\{0,+,-\}$
(see Example \ref{EX:RTW-SU2}), via
$\triangle_-=\D_{12}$, $\triangle_+=\D_{21}$, 
and $\triangle_0=\D_{11}-\D_{22}$, 
while $\D_{11}+\D_{22} \in \diff^2(\widehat G)$, 
and satisfy the Leibniz rule \eqref{EQ:G-leinbiz-1}.
Furthermore, the elementary (homogeneous) first 
order differential operators 
$\partial_0,\partial_\pm\in\Psi^1(\G)$ 
and the Laplacian
$\mathcal L=-\partial_0^2-
\frac12(\partial_-\partial_++\partial_+\partial_-)\in 
\Psi^2(\G)$ satisfy
\begin{center}
\begin{tabular}{c|c|c|c|c}
\raisebox{3ex}[3ex][1.5ex]{~}$ $ ~&~ $ \sigma_{\partial_0} 
$ ~&~ $ \sigma_{\partial_+} $ ~&~ $ \sigma_{\partial_-} 
$ ~&~ $ \sigma_{\mathcal L} $
 \\\hline
\raisebox{3ex}[3ex][1.5ex]{~} $\D_{11} $ ~&~ $  \frac12 \sigma_I  $ ~&~ $  0 $ ~&~ $ 0 $ ~&~ $ -\sigma_{\partial_0}+\frac14  \sigma_I$ \\\hline
\raisebox{3ex}[3ex][1.5ex]{~}$\D_{12} $ ~&~ $0$ ~&~ $\sigma_I$ ~&~ $0$ ~&~ $-\sigma_{\partial_-}$ \\\hline
\raisebox{3ex}[3ex][1.5ex]{~}$\D_{21} $ ~&~$0$ ~&~ $0$ ~&~ $\sigma_I$ ~&~ $-\sigma_{\partial_+}$ \\\hline
\raisebox{3ex}[3ex][1.5ex]{~}$\D_{22} $ ~&~$  -\frac12 \sigma_I  $ ~&~ $  0 $ ~&~ $ 0 $ ~&~ $ \sigma_{\partial_0}+\frac14 \sigma_I $ 
\end{tabular}
\end{center}
Invariant vector fields corresponding to the 
basis Pauli matrices in $\mathfrak{su}(2)$ are expressible as
$\mathrm D_1 = -\frac{\mathrm i}2(\partial_-+\partial_+) $, 
$\mathrm D_2=\frac12 (\partial_--\partial_+)$ and 
$\mathrm D_3= -\mathrm i \partial_0$. We refer to
\cite{RT-book,RT-groups} for the details on these operators and
for the explicit formulae for their symbols.

\begin{ex}
We consider $\mathrm D_3 + c$ for constants $c\in\C$. 
The corresponding symbol is given
by 
$$
   ( \sigma_{\mathrm D_3+c})_{mn}^\ell = [c - \mathrm i m]
   \delta_{mn}^\ell
$$
where, following \cite{RT-book}, 
$m,n\in\{-\ell,-\ell+1,\ldots,\ell-1,\ell\}$ run through half-integers
or integers, depending on whether $\ell$ is a half-integer or
an integer, respectively. 
We also use the notation $\delta_{mn}^\ell$ for the Kronecker
delta, i.e. $\delta_{mn}^\ell=1$ for $m=n$ and
$\delta_{mn}^\ell=0$ for $m\not=n$.
This symbol $\sigma_{\mathrm D_3+c}$
is invertible for all $\ell$ if (and only if) $c$ 
is not an imaginary half-integer, 
$\mathrm i c \not\in\frac12\mathbb Z$. The inverse satisfies 
$(\mathrm D_3+c)^{-1} \in \Op(\mathscr S^0_{0,0}(\G))$, based on
\begin{align*}
  &\D_{11} \sigma_{(\mathrm D_3+c)^{-1}} = 
  - \textstyle\frac12  \sigma_{(\mathrm D_3+c)^{-1}} 
  \sigma_{(\mathrm D_3+c+\frac12)^{-1}}, \\
   &\D_{22} \sigma_{(\mathrm D_3+c)^{-1}} = 
   \textstyle\frac12  \sigma_{(\mathrm D_3+c)^{-1}} 
   \sigma_{(\mathrm D_3+c-\frac12)^{-1}},
\end{align*}
and $\D_{12}\sigma_{(\mathrm D_3+c)^{-1}} =
\D_{21}\sigma_{(\mathrm D_3+c)^{-1}} =0,$ 
combined with Leibniz' rule. As a consequence we see 
that the operator $\mathrm D_3+c$ satisfies sub-elliptic 
estimates with loss of one derivative and are thus globally hypoelliptic.

The statement is sharp: the spectrum of $\mathrm D_3$ 
consists of all imaginary half-integers and all 
eigenspaces are infinite-dimensional 
(stemming from the fact that each such 
imaginary half-integer hits infinitely many 
representations for which $-\ell\le \mathrm i c\le \ell$ 
and $\mathrm i c+\ell\in\mathbb Z$). In particular, 
eigenfunctions can be irregular, e.g., distributions 
which do not belong to certain negative order Sobolev spaces.
\end{ex}

\begin{ex}\label{EX:RTW-ex53}
By conjugation formulae in \cite{RT-groups,RT-book} 
this example can be extended to give examples of more general
left-invariant operators. Namely, let $X\in\mathfrak g$ and let
$\partial_X$ be the derivative with respect to the vector
field $X$. Then the operator 
$\partial_X+c$ is globally hypoelliptic
if and only if $\mathrm i c \not\in\frac12\mathbb Z$.

By conjugation again, we can readily see that in a suitably
chosen basis in the representation spaces, the symbol of
the operator $\partial_X^k$, $k\in\N$, is
$$
   ( \sigma_{\mathrm \partial_X^k})_{mn}^\ell = 
   [(- \mathrm i m)^k]
   \delta_{mn}^\ell.
$$
Denoting $(\frac12\Z)^k=\{2^{-k}l^k: l\in\Z\}$, we see that
the operator $\partial_X^k+c$ is globally hypoelliptic if
and only if $c\not\in -(-\mathrm i)^k (\frac12\Z)^k$.
In particular, for $k=2$, the second order differential
operator $\partial_X^2+c$ is globally hypoelliptic if
and only if $4c\not\in \{l^2: l\in\Z\}$. 

It is also interesting
to note that the operators $\partial_X^2+a(x)\partial_X$
are not globally hypoelliptic for all complex-valued
functions $a\in C^\infty(\SU2)$. Indeed, in a suitably chosen
basis of the representation spaces in
$\widehat{\SU2}$ their symbols are diagonal, and
$(\sigma_{\partial_X^2+a(x)\partial_X})^\ell_{mn}=0$ for
$\ell\in\N$ and $m=n=0$. 
For simplicity in the following construction we assume
that $\partial_X=\mathrm D_3$.
If we set
$\widehat{f}(\ell)_{mn}=1$ for $\ell\in\N$ and
$m=n=0$, and to be zero otherwise, we have
$(\sigma_{\partial_X^2+a(x)\partial_X})\widehat{f}=0$,
so that $(\partial_X^2+a(x)\partial_X)f=0$ by
\eqref{EQ:RTW-quant}. By the Fourier inversion formula
we see that $\widehat{f}$ is the Fourier transform of the
function
$$
 f(x)=\sum_{\ell\in\N_0} (2\ell+1) t^\ell_{00}(x),
$$
so that $f\in\mathcal D'(\SU2)$ is not smooth but
$(\partial_X^2+a(x)\partial_X)f=0$.

We note that operators from Example \ref{EX:RTW-ex53}
are not covered by 
H\"ormander's sum of the squares theorem
\cite{Ho1}. Although the
following example of the sub-Laplacian is covered by 
H\"ormander's theorem,
Theorem \ref{THM:RTW-hypoellipticity} provides its
inverse as a pseudo-differential operator with
the global matrix-valued symbol in the class
$\mathscr S^{-1}_{\frac12,0}(G)$.
\end{ex}

\begin{ex}
{\sl Sub-Laplacian.} The sub-Laplacian on $\SU2$ is 
defined as 
$\mathcal L_s = \mathrm D_1^2 + \mathrm D_2^2\in\Psi^2(\G)$
and has the symbol
$$
   (\sigma_{\mathcal L_s})_{mn}^\ell  = 
   [m^2 - \ell (\ell+1)] \delta_{mn}^\ell.
$$
This symbol is invertible for all $\ell>0$. 
It satisfies the hypoellipticity assumption 
of Theorem~\ref{THM:RTW-hypoellipticity} with $\rho=\frac12$. 
For this we note that 
$\D_{11} \sigma_{\mathcal L_s} = \D_{22}\sigma_{\mathcal L_s} = 0$ 
and $\D_{12} \sigma_{\mathcal L_s} = -\sigma_{\partial_-}$, 
$\D_{21}\sigma_{\mathcal L_s} = - \sigma_{\partial_+}$. 
Estimate \eqref{EQ:RTW-he2} is evident for $\alpha=0$ 
and $|\alpha|\ge 2$. Hence it remains to check that
$$
 \|( \sigma_{\mathcal L_s})^{-1} \D_{12} 
 \sigma_{\mathcal L_s}\| = 
 \|(\sigma_{\mathcal L_s})^{-1} \sigma_{\partial_-}\|
 \lesssim \ell^{-1/2}
$$
(and the analogous statement for $\D_{21}$). This can be done by direct computation,
 $$
 ( \sigma_{\mathcal L_s}^{-1} \D_{12} \sigma_{\mathcal L_s})^\ell_{mn}
 = \left[\frac{\sqrt{(\ell+n)(\ell-n+1)}}{(n-1)^2-\ell(\ell+1)}\right]\delta^\ell_{m,n-1},
 $$
the entries corresponding to $\ell=-n$ vanish, 
while all others are bounded by
$\sqrt {2 /\ell}$. Hence, $\mathcal L_s$ has 
a pseudo-differential parametrix 
$\mathcal L_s^\sharp \in \Op(\mathscr S^{-1}_{\frac12,0}(\G)$).
In particular, 
we automatically get the subelliptic estimate
$\|f\|_{H^{s}(G)}\leq
C\|\mathcal L_sf\|_{H^{s-1}(G)}$ for all $s\in\R$.
\end{ex}

\begin{ex} {\sl Heat operator.}
We consider the analogue of the heat operator on $\SU2$,
namely the operator $H=\mathrm D_3-\mathrm D_1^2-
\mathrm D_2^2$. Writing $H=\mathrm D_3-\Lap_s$, we see that
its symbol is
$$
   (\sigma_{H})_{mn}^\ell  = 
   [-\mathrm i m - m^2 + \ell (\ell+1)] \delta_{mn}^\ell.
$$
Consequently, using the differences for the sub-Laplacian,
we see that $\D_{11}\sigma_H=-\frac{\mathrm i}{2}\sigma_I$,
$\D_{22}\sigma_H=\frac{\mathrm i}{2}\sigma_I$,
$\D_{12}\sigma_H=-\D_{12}\sigma_{\Lap_s}=\sigma_{\partial_-}$,
$\D_{21}\sigma_H=-\D_{21}\sigma_{\Lap_s}=\sigma_{\partial_+}$.
By an argument similar to that for the sub-Laplacian, the
operator $H$ has a parametrix 
$H^\sharp\in\Op(\mathscr S^{-1}_{\frac12,0}(\G))$.
\end{ex}

\begin{ex} {\sl Schr\"odinger operator.}
We consider the analogue of the Schr\"odinger operator on $\SU2$,
namely the operators $S_\pm=\pm\mathrm i \mathrm D_3-\mathrm D_1^2-
\mathrm D_2^2$. We will treat both $\pm$ cases simultaneously
by keeping the sign $\pm$ or $\mp$. Writing 
$S_\pm=\pm\mathrm i \mathrm D_3-\Lap_s$, we see that
their symbols are
$$
   (\sigma_{S_\pm})_{mn}^\ell  = 
   [\pm m - m^2 + \ell (\ell+1)] \delta_{mn}^\ell.
$$
Consequently, we see that these symbols are not
invertible. Analogously to Example \ref{EX:RTW-ex53}
we can construct distributions in the null spaces of $S_\pm$.
Namely, if we define $\widehat{f_\pm}(\ell)_{mn}=1$
for $m=n=\mp\ell$, and to be zero otherwise, then we
readily see that the functions
$f_\pm(x)=\sum_{\ell\in\frac12\N_0} (2\ell+1)
t^\ell_{\mp\ell \mp\ell}(x)$
satisfy $S_\pm f_\pm=0$, but $f_\pm\in\mathcal D'(G)$
are not smooth.

On the other hand, the operators $S_\pm+c$ become globally
hypoelliptic
for almost all constants $c$. Let us consider the case $S_-$, the
argument for $S_+$ is similar. Since the symbol of
$S_-+c$ is 
$$
   (\sigma_{S_-+c})_{mn}^\ell  = 
   [c- m - m^2 + \ell (\ell+1)] \delta_{mn}^\ell,
$$
we see that it is invertible if and only if
$c\not\in\{\ell(\ell+1)-m(m+1): \ell\in\frac12\N,
\ell-m\in\Z, |m|\leq \ell\}.$ If, for example,
$c>0$ or $c\in\C\backslash\R$, then 
$\|\sigma_{(S_-+c)^{-1}}\|_{op}\leq C$, and hence
$S_-+c$ is globally hypoelliptic.
\end{ex}

\begin{ex} {\sl D'Alembertian.}
We consider $W=\mathrm D_3^2-\mathrm D_1^2-\mathrm D_2^2$.
Writing $W=\mathrm D_3^2-\mathcal L_s=2\mathrm D_3^2-\Lap$
with the Laplace operator $\Lap=
\mathrm D_1^2+\mathrm D_2^2+\mathrm D_3^2$, 
we see that the symbol of
$W\in\Psi^2(G)$ is given by 
$$
   (\sigma_{W})_{mn}^\ell  = 
   [-2m^2 + \ell (\ell+1)] \delta_{mn}^\ell.
$$
If $\ell$ is a half-integer, we 
easily see that $|-2m^2+\ell(\ell+1)|\geq \frac14$ for all 
$\ell-m\in\Z$, $|m|\leq \ell$. However, if $\ell$ is an integer,
the triangular numbers $\frac{\ell(\ell+1)}{2}$ are squares if
$\ell = \ell_k =  \lfloor(3+2\sqrt{2})^k / 4\rfloor$, 
$k\in\N$, see \cite[A001109]{Sloane}. Denoting by $m_k$ the corresponding solution to
$2m_k^2=\ell_k(\ell_k+1)$, we see that the function
$$
 f(x)=\sum_{k=1}^\infty (2\ell_k+1) t^\ell_{m_k m_k}(x)
$$
satisfies $Wf=0$ but $f\in \mathcal D'(G)$ is not smooth.

Since the only non-zero triangular number which is also a cube is $1$ (see, e.g., \cite{Wolfram:}) we see that the operator
$2\mathrm i \mathrm D_3^3-\Lap$ is globally hypoelliptic.
\end{ex}

\begin{ex} 
We consider $P=\mathrm D_1^2 - \mathrm D_2^2 = 
- \frac12 (\partial_-^2 + \partial_+^2)$. 
The corresponding symbol satisfies
\begin{align*}
2 (\sigma_P)^\ell_{mn} =&\sqrt{(\ell+n)(\ell+n-1)
(\ell-n+2)(\ell-n+1)} \delta^\ell_{m,n-2} \\
&+ \sqrt{(\ell+m)(\ell+m-1)(\ell-m+2)(\ell-m+1)} 
\delta^\ell_{m-2,n}.
\end{align*}
Furthermore, by Leibniz' rule we deduce 
$\D_{11}\sigma_P = \D_{22}\sigma_P=0$ and 
$\D_{12} \sigma_P = - \sigma_{\partial_+}$, 
$\D_{21}\sigma_P=-\sigma_{\partial_-}$.
The symbol $\sigma_P$ is invertible only for 
$\ell+\frac12\in 2\mathbb N_0$. 
Therefore, $P$ is neither elliptic 
nor does it satisfy subelliptic estimates.
\end{ex}


\end{document}